\documentclass[final,1p,times]{elsarticle}

\usepackage{amssymb}
 \usepackage{amsthm}
\usepackage{amscd}
\usepackage{amsmath}
\usepackage{amsfonts}
\usepackage{amssymb}
\usepackage{graphicx}
\newtheorem{theorem}{Theorem}
\newtheorem{proposition}[theorem]{Proposition}
\newtheorem{lemma}[theorem]{Lemma}
\newtheorem{rmk}[theorem]{Remark}
\newtheorem{corollary}[theorem]{Corollary}
\usepackage{mathrsfs}
\usepackage{titletoc}

\newcommand{\la}{\Delta}
\def\na{\nabla}

\renewcommand{\d}{\delta}

\newcommand{\ra}{\rightarrow}
\newcommand{\p}{\partial}
\newcommand{\f}{\frac}
\def\a{\alpha}
\def\lam{\lambda}

\def\e{\epsilon}

\def\S{\Sigma}

\newcommand{\be}{\begin{equation}}
\renewcommand{\ra}{\rightarrow}
\newcommand{\ee}{\end{equation}}
\newcommand{\bea}{\begin{eqnarray}}
\newcommand{\eea}{\end{eqnarray}}
\newcommand{\bna}{\begin{eqnarray*}}
\newcommand{\ena}{\end{eqnarray*}}

\renewcommand{\le}{\left}
\newcommand{\ri}{\right}

\journal{***}

\begin{document}

\begin{frontmatter}

\title{Existence of solutions to a class of Kazdan-Warner equations on compact Riemannian surface}

\author{Yunyan Yang}
 \ead{yunyanyang@ruc.edu.cn}
\author{Xiaobao Zhu}
\ead{zhuxiaobao@ruc.edu.cn}
\address{ Department of Mathematics,
Renmin University of China, Beijing 100872, P. R. China}

\begin{abstract}
Let $(\Sigma,g)$ be a compact Riemannian surface without boundary and $\lambda_1(\Sigma)$ be the first eigenvalue of the Laplace-Beltrami
operator $\Delta_g$. Let $h$ be a positive smooth function on $\Sigma$.
Define a functional
$$J_{\alpha,\beta}(u)=\f{1}{2}\int_\Sigma(|\nabla_gu|^2-\alpha u^2)dv_g-\beta\log\int_\Sigma he^udv_g$$
on a function space $\mathcal{H}=\le\{u\in W^{1,2}(\Sigma): \int_\Sigma udv_g=0\ri\}$. If $\alpha<\lambda_1(\Sigma)$ and $J_{\alpha,8\pi}$ has no minimizer
on $\mathcal{H}$, then we calculate the infimum of $J_{\alpha,8\pi}$ on $\mathcal{H}$ by using the method of blow-up analysis.
As a consequence, we give a sufficient condition under which a Kazdan-Warner equation has a solution.
If $\alpha\geq \lambda_1(\Sigma)$, then $\inf_{u\in\mathcal{H}}J_{\alpha,8\pi}(u)=-\infty$. If $\beta>8\pi$, then for any
$\alpha\in\mathbb{R}$, there holds $\inf_{u\in\mathcal{H}}J_{\alpha,\beta}(u)=-\infty$.
Moreover, we consider the same problem in the case that $\alpha$ is large, where higher order eigenvalues are involved.
\end{abstract}


\end{frontmatter}

\titlecontents{section}[0mm]
                       {\vspace{.2\baselineskip}}
                       {\thecontentslabel~\hspace{.5em}}
                        {}
                        {\dotfill\contentspage[{\makebox[0pt][r]{\thecontentspage}}]}
\titlecontents{subsection}[3mm]
                       {\vspace{.2\baselineskip}}
                       {\thecontentslabel~\hspace{.5em}}
                        {}
                       {\dotfill\contentspage[{\makebox[0pt][r]{\thecontentspage}}]}

\setcounter{tocdepth}{2}

\section{Introduction and main results}
Let $(\Sigma,g)$ be a compact Riemannian surface without boundary, $W^{1,2}(\Sigma)$ be the usual Sobolev space.
Define a function space
\be\label{H}\mathcal{H}=\le\{u\in W^{1,2}(\Sigma):\,\int_\Sigma udv_g=0\ri\}.\ee
Let $h$ be a positive smooth function on $\Sigma$ and $J_{\beta}: W^{1,2}(\Sigma)\ra\mathbb{R}$ be a functional defined by
\be\label{2}J_{\beta}(u)=\f{1}{2}\int_\Sigma|\nabla_gu|^2dv_g-\beta\log\int_\Sigma he^udv_g,\ee
where $\nabla_gu$ denotes the gradient of $u$ and $dv_g$ denotes the volume element of $(\Sigma,g)$.
In view of the Trudinger-Moser inequality due to Fontana \cite{Fontana},
$J_\beta$ has a minimizer on $\mathcal{H}$ for any $\beta<8\pi$; while in the case $\beta=8\pi$, the situation becomes subtle.
Using a method of blow-up analysis, Ding-Jost-Li-Wang \cite{DJLW} proved that if $J_{8\pi}$ has no minimizer on $\mathcal{H}$, then
\be\label{inf-DJLW}\inf_{u\in\mathcal{H}}J_{8\pi}(u)=-8\pi-8\pi\log\pi-4\pi
\max_{p\in\Sigma}(A_p+2\log h(p)),\ee
where $A_p=\lim_{x\ra p}(G_p(x)+4\log r)$ is a constant, $r$ denotes the geodesic distance between $x$ and $p$, $G_p$ is a Green function
satisfying
$$
\begin{cases}
      \Delta_g G_p = 8\pi\d_{p}-\f{8\pi}{{\rm Vol}_g(\Sigma)} \\[1.2ex]
      \int_{\Sigma}G_p dv_{g}=0,
\end{cases}
$$
and $\Delta_g$ is the Laplace-Beltrami operator.
Moreover, they give a geometric condition under which $J_{8\pi}$ has a minimizer on $\mathcal{H}$.
Clearly the minimizer is a solution of a Kazdan-Warner equation \cite{Kazdan-Warner}, namely
\be\label{KZ-0}\Delta_g u=\f{8\pi he^u}{\int_\Sigma he^udv_g}-\f{8\pi}{{\rm Vol}_g(\Sigma)}.\ee

 We let $\lambda_1(\Sigma)$ be the first eigenvalue of $\Delta_g$, say
\be\label{eig}\lambda_1(\Sigma)=\inf_{u\in \mathcal{H},\,\int_\Sigma u^2dv_g=1}
\int_\Sigma|\nabla_gu|^2dv_g.\ee
It follows from the Poincar\'e inequality that if $\alpha<\lambda_1(\Sigma)$, then
$$\|u\|_{1,\alpha}=\le(\int_\Sigma(|\nabla_gu|^2-\alpha u^2)dv_g\ri)^{1/2}$$
defines a Sobolev norm on $\mathcal{H}$. In a previous work \cite{YangJDE}, using the method of blow-up analysis,
we proved the following:  for any $\alpha<\lambda_1(\Sigma)$,
 there holds
\be\label{strong-TM}\sup_{u\in \mathcal{H},\,\|u\|_{1,\alpha}\leq 1}\int_\Sigma e^{4\pi u^2}dv_g<\infty\ee
and the supremum is attained. As a consequence of (\ref{strong-TM}),
there exists some constant $C$ depending only on $(\Sigma,g)$ and $\alpha<\lambda_1(\Sigma)$ such that for all $u\in\mathcal{H}$,
\be\label{8}\f{1}{2}\|u\|_{1,\alpha}^2-8\pi \log\int_\Sigma he^udv_g\geq -C.\ee
This improves the Trudinger-Moser inequality of the weak form, namely (\ref{8}) in the case $\alpha=0$.
We refer
the reader to \cite{A-D,YangTran,Lu-Yang,doo1,doo2,Tintarev,YangJGA,YangZhu,YangZhu2}
for related works involving the norms $\|u\|_{1,\alpha}$.

Our aim in this paper is to achieve an analog of (\ref{inf-DJLW}).
More precisely, we consider functionals
\be\label{Jab}J_{\alpha,\beta}(u)=\f{1}{2}\int_\Sigma(|\nabla_gu|^2-\alpha u^2)dv_g-\beta\log\int_\Sigma he^udv_g.\ee
Obviously, when $\alpha=0$, $J_{\alpha,\beta}$ reduces to $J_\beta$ defined as in (\ref{2}).
Our first result reads
\begin{theorem}\label{main-theorem}
Let $(\Sigma,g)$ be a compact Riemannian surface without boundary, $h$ be a positive smooth function on $\Sigma$,
and $\mathcal{H}$, $\lambda_1(\Sigma)$ and $J_{\alpha,\beta}$ be defined
as in (\ref{H}), (\ref{eig}) and (\ref{Jab}) respectively. Then we have the following three assertions:\\
$(i)$ If $\alpha<\lambda_1(\Sigma)$ and $J_{\alpha,8\pi}$
has no minimizer in $\mathcal{H}$, then there holds
$$\inf_{u\in\mathcal{H}}J_{\alpha,8\pi}(u)=-8\pi-8\pi\log\pi-4\pi
\max_{p\in\Sigma}(A_{p}+2\log h(p)),$$
where $A_{p}=\lim_{x\ra p}(G_{p}(x)+4\log r)$ is a constant, $r$ denotes the geodesic distance between $x$ and $p$,
$G_{p}$ is a Green function
satisfying
\be\label{green-alpha}
\begin{cases}
       \Delta_g G_{p}-\alpha G_{p}= 8\pi\d_{p}-\f{8\pi}{{\rm Vol}_g(\Sigma)} \\[1.2ex]
       \int_{\Sigma}G_{p} dv_{g}=0;
\end{cases}
\ee
$(ii)$ If $\alpha\geq \lambda_1(\Sigma)$, then $\inf_{u\in\mathcal{H}}J_{\alpha,8\pi}(u)=-\infty;$\\
$(iii)$ If $\beta>8\pi$, then for any $\alpha\in\mathbb{R}$, we have
$\inf_{u\in\mathcal{H}}J_{\alpha,\beta}(u)=-\infty$.
\end{theorem}
Since the Euler-Lagrange equation of a minimum point of $J_{\alpha,8\pi}$ on $\mathcal{H}$ is
\be\label{Euler}\Delta_g u-\alpha u=\f{8\pi he^u}{\int_\Sigma he^udv_g}-\f{8\pi}{{\rm Vol}_g(\Sigma)},\ee
 an application of $(i)$ of Theorem \ref{main-theorem} is the following:
 \begin{corollary}\label{Cor-1}
 For any $\alpha<\lambda_1(\Sigma)$,
 if
 \be\label{c-1}\inf_{u\in\mathcal{H}}J_{\alpha,8\pi}(u)\not=-8\pi-8\pi\log\pi-4\pi
\max_{p\in\Sigma}(A_{p}+2\log h(p)),\ee
then the Kazdan-Warner equation (\ref{Euler}) has a solution $u\in\mathcal{H}$.
 \end{corollary}

 As in \cite{YangJDE}, we consider the case that $\alpha$
  is allowed to be larger than $\lambda_1(\Sigma)$. Precisely, we let
  $\lambda_1(\Sigma)<\lambda_2(\Sigma)<\cdots$ be all distinct eigenvalues
  of $\Delta_g$, $E_{\lambda_\ell(\Sigma)}$
 be the eigenfunction space with respect to $\lambda_\ell(\Sigma)$, namely
 \be\label{eigenfunct}E_{\lambda_\ell(\Sigma)}=\le\{u\in\mathcal{H}:\,\Delta_gu=\lambda_\ell(\Sigma)u\ri\},\ee
 and
 $E_\ell=E_{\lambda_1(\Sigma)}\oplus E_{\lambda_2(\Sigma)}\oplus \cdots\oplus E_{\lambda_\ell(\Sigma)}$,
 $\ell=1,2,\cdots$. Define
 \be\label{E-lp}E_\ell^\perp=\le\{u\in \mathcal{H}: \int_\Sigma uvdv_g=0,\,\forall v\in  E_\ell\ri\}.\ee
 Now we state an analog of Theorem \ref{main-theorem} as follows:
 \begin{theorem}\label{main-theorem-2}
Let $(\Sigma,g)$, $h$, $\mathcal{H}$ and $J_{\alpha,\beta}$ be as in Theorem \ref{main-theorem}, $\lambda_\ell(\Sigma)$
be the $\ell$-th eigenvalue of the Laplace-Beltrami operator, $E_\ell$ and
$E_\ell^\perp$ be defined as in (\ref{E-lp}).
Then we have the following three assertions:\\
$(i)$ If $\alpha<\lambda_{\ell+1}(\Sigma)$ and $J_{\alpha,8\pi}$
has no minimizer in $E_\ell^\perp$, then there holds
$$\inf_{u\in E_\ell^\perp}J_{\alpha,8\pi}(u)=-8\pi-8\pi\log\pi-4\pi
\max_{p\in\Sigma}(A_{\alpha,p}+2\log h(p)),$$
where $A_{\alpha,p}=\lim_{x\ra p}(G_{\alpha,p}(x)+4\log r)$ is a constant, $r$ denotes the geodesic distance between $x$ and $p$,
$G_{\alpha,p}$ is a Green function
satisfying
\be\label{green-alpha}
\begin{cases}
       \Delta_g G_{\alpha,p}-\alpha G_{\alpha,p}=8\pi\d_{p}-\f{8\pi}{{\rm Vol}_g(\Sigma)} \\[1.2ex]
       \int_{\Sigma}G_{\alpha,p}v dv_{g}=0,\,\forall v\in E_\ell;
\end{cases}
\ee
$(ii)$ If $\alpha\geq \lambda_{\ell+1}(\Sigma)$, then
$\inf_{u\in E_\ell^\perp}J_{\alpha,8\pi}(u)=-\infty;$\\
$(iii)$ If $\beta>8\pi$, then for any $\alpha\in\mathbb{R}$, we have
$\inf_{u\in E_\ell^\perp}J_{\alpha,\beta}(u)=-\infty$.
\end{theorem}

Similar to Corollary \ref{Cor-1}, we have the following:
 \begin{corollary}\label{Cor-2}
 For any $\alpha<\lambda_{\ell+1}(\Sigma)$,
 if
 \be\label{c-2}\inf_{u\in E_\ell^\perp}J_{\alpha,8\pi}(u)\not=-8\pi-8\pi\log\pi-4\pi
\max_{p\in\Sigma}(A_{\alpha,p}+2\log h(p)),\ee
then the Kazdan-Warner equation (\ref{Euler}) has a solution $u\in E_\ell^\perp$.
 \end{corollary}

 We remark that any geometric hypothesis under which (\ref{c-1}) or (\ref{c-2}) holds would be extremely interesting.
 When $\alpha=0$, a geometric condition was given by Ding-Jost-Li-Wang \cite{DJLW}. Generally it is difficult to be obtained
 possibly because the Green function has at most $C^{1,\gamma}$-regularity for some $0<\gamma<1$ in presence of $\alpha$.

 Now we describe our method. For the proof of $(ii)$ and $(iii)$ of Theorems \ref{main-theorem} and \ref{main-theorem-2},
 we shall construct suitable function sequences. To prove $(i)$ of Theorems \ref{main-theorem} and \ref{main-theorem-2},
 we use the blow-up scheme proposed by Ding-Jost-Li-Wang \cite{DJLW}.
 Our analysis is different from that of Ding-Jost-Li-Wang \cite{DJLW} at least in three points: Let $u_\epsilon$
 be a minimizer of $J_{\alpha,8\pi(1-\epsilon)}$ and $u_\epsilon(x_\epsilon)=\max_{\Sigma}u_\epsilon\ra+\infty$.
 One is that before understanding the exact asymptotic behavior of $u_\epsilon$ near the blow-up point,
 we must prove $u_\epsilon(x)/u_\epsilon(x_\epsilon)= 1+o_\epsilon(1)$ in $B_{Rr_\epsilon}(x_\epsilon)$, where $r_\epsilon\ra 0$ is an appropriate
 sequence of positive numbers, $R>0$ is fixed, and $o_\epsilon(1)\ra 0$ uniformly in $B_{Rr_\epsilon}(x_\epsilon)$;
 The other is that in the process of deriving lower bound of $J_{\alpha,8\pi}$,
 we estimate the energy $\int_\Sigma |\nabla_gu_\epsilon|^2dv_g$ on two regions $B_{Rr_\epsilon}$ and
 $\Sigma\setminus B_{Rr_\epsilon}$
   instead of three regions $B_{Rr_\epsilon}$, $B_\delta\setminus B_{Rr_\epsilon}$ and
   $\Sigma\setminus B_{\delta}$, which simplifies the calculation in \cite{DJLW}; The third is in the final step (test function
 computation), we construct a sequence of test functions different from that of \cite{DJLW}.

  Before ending this introduction, we mention several related works also based on the blow-up scheme in \cite{DJLW}.
  Ni \cite{Ni} considered the mean field equation with critical parameter in a planar domain.
  Zhou \cite{Zhou} obtained existence of solution to the mean field equation for the equilibrium turbulence.
  Liu-Wang \cite{LiuWang} studied the equation (\ref{KZ-0}) with an extra
    drifting term $\nabla\phi\cdot\nabla u$. Mancini \cite{Mancini} proved an Onofri inequality.

 Throughout this paper, $o_\epsilon(1)\ra 0$ as $\epsilon\ra 0$, $o_R(1)\ra 0$ as $R\ra\infty$, and so on.
 We do not distinguish sequence and subsequence and often denote various constants by the same $C$.
   The remaining part of this paper is organized as follows: In Section \ref{sec1}, we prove Theorem \ref{main-theorem};
    In Section \ref{sec2}, we give the proof of Theorem \ref{main-theorem-2}.

\section{Proof of Theorem \ref{main-theorem}}\label{sec1}

 The proof of $(ii)$ and $(iii)$ is easy and will be shown first.

\subsection{Proof of $(ii)$ of Theorem \ref{main-theorem}}

Let $\alpha\geq \lambda_1(\Sigma)$ be fixed and $E_{\lambda_1(\Sigma)}$ be the eigenfunction space defined as in (\ref{eigenfunct}).
Take $0\not\equiv u_0\in E_{\lambda_1(\Sigma)}$. Obviously we have
$$\int_\Sigma|\nabla_gu_0|^2dv_g=\lambda_1(\Sigma)\int_\Sigma u_0^2dv_g$$
and thus
\be\label{less0}\int_\Sigma(|\nabla_gu_0|^2-\alpha u_0^2)dv_g\leq 0.\ee
Since $\int_\Sigma u_0dv_g=0$, there exists $x_0\in \Sigma$ and $\delta>0$ such that
\be\label{geq0} u_0(x)\geq u_0(x_0)/2>0,\quad\forall
x\in B_\delta(x_0).\ee
It follows from (\ref{less0}) and (\ref{geq0}) that for any $t>0$,
\bna
J_{\alpha,8\pi}(t u_0)&=&\f{t^2}{2}\int_\Sigma(|\nabla_gu_0|^2-\alpha u_0^2)dv_g-8\pi\log\int_\Sigma
he^{tu_0}dv_g\\
&\leq&-8\pi\log\int_{B_\delta(x_0)}he^{tu_0}dv_g\\
&\leq&-4\pi tu_0(x_0)-8\pi\log\int_{B_\delta(x_0)}hdv_g.
\ena
 Hence $J_{\alpha,8\pi}(tu_0)\ra-\infty$ as $t\ra+\infty$ and the desired result follows immediately. $\hfill\Box$

\subsection{Proof of $(iii)$ of Theorem \ref{main-theorem}}

 Let $\beta>8\pi$ be fixed and $i_g(\Sigma)$ be the injectivity radius of $(\Sigma,g)$.
 Fix some point $p\in\Sigma$. Let $r$, $0<r<i_g(\Sigma)/2$, be a real number to be determined later. Take a sequence of
 functions
$$M_k=M_k(x,r)=\le\{
\begin{array}{lll}
\log k&{\rm when} &\rho\leq rk^{-\f{1}{4}}\\[1.5ex]
4\log\f{r}{\rho} &{\rm when} &rk^{-\f{1}{4}}<\rho\leq r\\[1.5ex]
0&{\rm when}& \rho>r,
\end{array}
\ri.$$
where $x\in\Sigma$, $\rho$ denotes the geodesic distance between $x$ and $p$, and $k=2,3\cdots$. One calculates
\bea\label{23}
&&\int_\Sigma|\nabla_g M_k|^2dv_g=(1+O(r))8\pi\log k,\label{g-1}\\
&&\int_\Sigma M_k^j d v_g=O(1),\quad j=1,2.\label{g-2}
\eea
Choose $\zeta\in C_0^1(\Sigma\setminus B_{2i_g(\Sigma)}(p))$ satisfying $\zeta\geq 0$  and
$\zeta\not\equiv 0$ on $\Sigma\setminus B_{2i_g(\Sigma)}(p)$.
We set
\be\label{M-tilde}\widetilde{M}_k=\widetilde{M}_k(x,r)=\le\{
\begin{array}{lll}
M_k(x,r)&{\rm for} &x\in B_{2i_g(\Sigma)}(p)\\[1.5ex]
t_k\zeta(x)&{\rm for}& x\in\Sigma\setminus B_{2i_g(\Sigma)}(p),
\end{array}
\ri.\ee
where $t_k\in \mathbb{R}$ is chosen so that $\widetilde{M}_k\in\mathcal{H}$. As a consequence,
$$\label{tk}t_k=-\f{1}{\int_\Sigma \zeta dv_g}\int_\Sigma M_k dv_g=O(1).$$
It follows from (\ref{23})-(\ref{M-tilde}) that
\be\label{e-1}\int_\Sigma|\nabla_g \widetilde{M}_k|^2dv_g-\alpha\int_\Sigma \widetilde{M}_k^2dv_g=(1+O(r))8\pi\log k+O(1)\ee
and
\be\label{e-2}\int_\Sigma he^{\widetilde{M}_k}dv_g\geq \int_{B_{rk^{-{1}/{4}}}(p)}he^{{M}_k}dv_g\geq
(\min_\Sigma h)\pi r^2 \sqrt{k}(1+o_k(1)).\ee
Combining (\ref{e-1}) and (\ref{e-2}), we have
$$J_{\alpha,\beta}(\widetilde{M}_k)\leq (1+O(r))4\pi\log k-\f{\beta}{2}\log k+O(1).$$
Note that $\beta>8\pi$. If $r$ is chosen sufficiently small, then we conclude
$$J_{\alpha,\beta}(\widetilde{M}_k)\ra-\infty \quad {\rm as}\quad k\ra+\infty.$$
This completes the proof of $(iii)$ of Theorem \ref{main-theorem}. $\hfill\Box$\\

In the remaining part of this section, we always assume $\alpha<\lambda_1(\Sigma)$. Since the proof of
$(i)$ of Theorem \ref{main-theorem} is very long, we sketch its outline as follows:
{\it Step 1}. For any $\epsilon$, $0<\epsilon<1$, there exists a minimizer
$u_\epsilon\in\mathcal{H}\cap C^1(\Sigma)$ for the subcritical functional $J_{\alpha,8\pi(1-\epsilon)}$.
By assumption that $J_{\alpha,8\pi}$ has no minimizer on $\mathcal{H}$,
we have $u_\epsilon(x_\epsilon)=\max_\Sigma u_\epsilon\ra+\infty$ as
$\epsilon\ra 0$. {\it Step 2}. We prove that $u_\epsilon(\exp_{x_\epsilon}(r_\epsilon x))-\max_\Sigma u_\epsilon\ra \varphi(x)$
in $C^1_{\rm loc}(\mathbb{R}^2)$, where $\exp_{x_\epsilon}$ denotes the exponential map on $(\Sigma,g)$, $r_\epsilon$ is an appropriate
scale, and $\varphi$ can be explicitly written out via a classification result of Chen-Li \cite{CL}.
Moreover, assuming $x_\epsilon\ra p$, we show that $u_\epsilon\ra G_p$ weakly in $W^{1,q}(\Sigma)$ for any $1<q<2$ and in
$C^1_{\rm{loc}}(\Sigma\setminus\{p\})$, where $G_p$ is a Green function on $(\Sigma,g)$. {\it Step 3}. Applying the maximum principle
to $u_\epsilon-G_{x_\epsilon}$ and using the
asymptotic behavior of $u_\epsilon$ derived in Step 2,
we obtain a lower bound of $J_{\alpha,8\pi}$ on $\mathcal{H}$. {\it Step 4}. We construct a sequence of functions $\phi_\epsilon\in \mathcal{H}$
such that $J_{\alpha,8\pi}(\phi_\epsilon)$ converges to the lower bound obtained in Step 3.

\subsection{Minimizers for subcritical functionals}\label{2.3}
We first prove that $\inf_{u\in\mathcal{H}}J_{\alpha,\beta}(u)$ is attained for any $\beta<8\pi$. Precisely we have
\begin{proposition}\label{Prop1}
For any $0<\epsilon<1$, there exists some function $u_\epsilon\in \mathcal{H}\cap C^1(\Sigma)$ such that
\be\label{sub}J_{\alpha,8\pi(1-\epsilon)}(u_\epsilon)=\inf_{u\in \mathcal{H}}J_{\alpha,8\pi(1-\epsilon)}(u).\ee
Moreover, $u_\epsilon$ satisfies the Euler-Lagrange equation
\be\label{EL}
\le\{\begin{array}{lll}\Delta_gu_\epsilon-\alpha u_\epsilon=8\pi(1-\epsilon)\le(\lambda_\epsilon^{-1}he^{u_\epsilon}-
\f{1}{{\rm Vol}_g(\S)}\ri)
\\[1.5ex]
\lambda_\epsilon=\int_\Sigma he^{u_\epsilon}dv_g,\quad \int_\Sigma u_\epsilon dv_g=0.\end{array}
\ri.\ee
\end{proposition}

\noindent{\it Proof.} Let $0<\epsilon<1$ be fixed. Take $u_j\in\mathcal{H}$ such that
\be\label{sequence}J_{\alpha,8\pi(1-\epsilon)}(u_j)\ra \inf_{u\in\mathcal{H}}J_{\alpha,8\pi(1-\epsilon)}(u)\ee
as $j\ra\infty$. Noting that
\be\label{bdd}\int_\Sigma he^{u_j}dv_g\leq \int_\Sigma he^{4\pi(1-\epsilon/2)\f{u_j^2}{\|u_j\|_{1,\alpha}^2}+
\f{\|u_j\|_{1,\alpha}^2}{16\pi(1-\epsilon/2)}}dv_g,\ee
we obtain by (\ref{sequence}) and (\ref{strong-TM}),
\bna
\inf_{u\in\mathcal{H}}J_{\alpha,8\pi(1-\epsilon)}(u)+o_j(1)&=&\f{1}{2}\|u_j\|_{1,\alpha}^2-8\pi(1-\epsilon)\log\int_\Sigma
he^{u_j}dv_g\\
&\geq&\f{1}{2}\|u_j\|_{1,\alpha}^2-\f{1-\epsilon}{2-\epsilon}\|u_j\|_{1,\alpha}^2-8\pi(1-\epsilon)\log
\int_\Sigma he^{4\pi(1-\epsilon/2)\f{u_j^2}{\|u_j\|_{1,\alpha}^2}}dv_g\\
&\geq&\f{\epsilon}{4}\|u_j\|_{1,\alpha}^2-C.
\ena
Hence $u_j$ is bounded in $\mathcal{H}$. We can assume without loss of generality that
$u_j$ converges to $u_\epsilon$ weakly in $\mathcal{H}$, strongly in $L^q(\Sigma)$ for any $q>0$ and
almost everywhere in $\Sigma$. Clearly
\be\label{weak-cont}\|u_\epsilon\|_{1,\alpha}^2\leq\lim_{j\ra\infty}\|u_j\|_{1,\alpha}^2.\ee
Moreover an analog of (\ref{bdd}) implies that $e^{|u_j|}$ is bounded in $L^p(\Sigma)$ for any $p>0$.
This together with the mean value theorem and the H\"older inequality,
\be\label{mean}\lim_{j\ra \infty}\int_\Sigma he^{u_j}dv_g=\int_\Sigma he^{u_\epsilon}dv_g.\ee
Combining (\ref{weak-cont}) and (\ref{mean}), we conclude (\ref{sub}).

Using a method of Lagrange multiplier, one easily gets (\ref{EL}), the Euler-Lagrange equation of
the minimizer $u_\epsilon$. Applying elliptic estimates to (\ref{EL}), we have  $u_\epsilon\in C^1(\Sigma)$.
$\hfill\Box$

\begin{lemma}\label{l>0}
$\liminf_{\epsilon\ra 0}\lambda_\epsilon>0$.
\end{lemma}
\noindent{\it Proof.} One may conclude the lemma by using the Jensen inequality. But we prefer a contradiction argument as below.
Clearly
\be\label{lll}J_{\alpha,8\pi(1-\epsilon)}(u_\epsilon)=\inf_{u\in \mathcal{H}}J_{\alpha,8\pi(1-\epsilon)}(u)\leq
J_{\alpha,8\pi(1-\epsilon)}(0)\leq 8\pi|\log\int_\Sigma hdv_g|.\ee
If $\liminf_{\epsilon\ra 0}\lambda_\epsilon=0$, then up to a subsequence
$$J_{\alpha,8\pi(1-\epsilon)}(u_\epsilon)=\f{1}{2}\|u_\epsilon\|_{1,\alpha}^2-8\pi(1-\epsilon)\log\lambda_\epsilon\ra +\infty,$$
which contradicts (\ref{lll}). Thus we get the desired result. $\hfill\Box$

\begin{lemma}\label{lim}
There holds
$$\lim_{\epsilon\ra 0}\inf_{u\in \mathcal{H}}J_{\alpha,8\pi(1-\epsilon)}(u)=\inf_{u\in\mathcal{H}}J_{\alpha,8\pi}(u).$$
\end{lemma}

\noindent{\it Proof.} Though the proof may be obvious for experts, we give the details here for reader's convenience.
On one hand, for any $\eta>0$, there exists some $u_\eta\in\mathcal{H}$
such that
$$J_{\alpha,8\pi}(u_\eta)<\inf_{u\in\mathcal{H}}J_{\alpha,8\pi}(u)+\eta.$$
Obviously we have
\bna J_{\alpha,8\pi}(u_\eta)=\lim_{\epsilon\ra 0}J_{\alpha,8\pi(1-\epsilon)}(u_\eta)
\geq\lim_{\epsilon\ra 0}J_{\alpha,8\pi(1-\epsilon)}(u_\epsilon)=\lim_{\epsilon\ra 0}
\inf_{u\in \mathcal{H}}J_{\alpha,8\pi(1-\epsilon)}(u).
\ena
Hence
$$\lim_{\epsilon\ra 0}
\inf_{u\in \mathcal{H}}J_{\alpha,8\pi(1-\epsilon)}(u)<\inf_{u\in\mathcal{H}}J_{\alpha,8\pi}(u)+\eta.$$
Since $\eta>0$ is arbitrary, we have
\be\label{leq}\lim_{\epsilon\ra 0}
\inf_{u\in \mathcal{H}}J_{\alpha,8\pi(1-\epsilon)}(u)\leq\inf_{u\in\mathcal{H}}J_{\alpha,8\pi}(u).\ee

On the other hand,
$$\inf_{u\in \mathcal{H}}J_{\alpha,8\pi}(u)\leq J_{\alpha,8\pi}(u_\epsilon)=\lim_{\nu\ra 0}J_{\alpha,8\pi(1-\nu)}(u_\epsilon).$$
Extracting diagonal sequence, we obtain
\be\label{geq}\inf_{u\in \mathcal{H}}J_{\alpha,8\pi}(u)\leq \lim_{\epsilon\ra 0}J_{\alpha,8\pi(1-\epsilon)}(u_\epsilon).\ee
Combining (\ref{leq}) and (\ref{geq}), we get the desired result. $\hfill\Box$

\begin{lemma}\label{lamd}
If $\lambda_\epsilon$ is bounded, then $u_\epsilon$ is bounded in $\mathcal{H}$ and a minimizer for
the functional $J_{\alpha, 8\pi}$ exists on the function space $\mathcal{H}$.
\end{lemma}

\noindent{\it Proof.}
Since $\lambda_\epsilon$ is a bounded sequence, it follows from (\ref{lll}) that
\bna
\f{1}{2}\|u_\epsilon\|_{1,\alpha}^2&=&J_{\alpha,8\pi(1-\epsilon)}(u_\epsilon)+8\pi(1-\epsilon)\log\int_\Sigma he^{u_\epsilon}dv_g\\
&\leq&8\pi|\log\int_\Sigma h d v_g|+8\pi(1-\epsilon)\log\lambda_\epsilon\\
&\leq& C.
\ena
Hence  $u_\epsilon$ is bounded in $\mathcal{H}$ and thus $e^{|u_\epsilon|}$ is bounded in $L^p(\Sigma)$ for any
$p>0$. Applying elliptic estimates to the equation (\ref{EL}), in view of Lemma \ref{l>0}, we have that $u_\epsilon$ converges to
some $u_0\in\mathcal{H}$ in $C^1(\Sigma)$. By Lemma \ref{lim}, we have
$$J_{\alpha,8\pi}(u_0)=\lim_{\epsilon\ra 0}J_{\alpha,8\pi(1-\epsilon)}(u_\epsilon)=\inf_{u\in\mathcal{H}}J_{\alpha,8\pi}(u).$$
Therefore $u_0$ is a minimizer of $J_{\alpha,8\pi}$. $\hfill\Box$\\

Denote
\be\label{maximum}c_\epsilon=\max_{x\in\Sigma}u_\epsilon(x)=u_\epsilon(x_\epsilon).\ee
\begin{proposition}\label{Prop2}
If $c_\epsilon$ is bounded from above, then $J_{\alpha,8\pi}$ has a minimizer in $\mathcal{H}$.
\end{proposition}
\noindent{\it Proof.} Multiplying both sides of the equation (\ref{EL}) by $u_\epsilon$, we have by using
Lemma \ref{l>0}, the assumption that $c_\epsilon$ is bounded from above and the Sobolev embedding theorem,
\bna
\|u_\epsilon\|_{1,\alpha}^2\leq C\int_\Sigma|u_\epsilon|dv_g\leq C\|u_\epsilon\|_{1,\alpha}.
\ena
This implies that $u_\epsilon$ is bounded in $\mathcal{H}$. Applying elliptic estimates to (\ref{EL}),
we conclude that $u_\epsilon$ converges to a minimizer of $J_{\alpha,8\pi}$ in $C^1(\Sigma)$. $\hfill\Box$

\subsection{Blow-up analysis}
We now analyze the asymptotic behavior of $u_\epsilon$.
By our assumption that $J_{\alpha,8\pi}$ has no minimizer on $\mathcal{H}$, in view of Lemma \ref{lamd}
and Proposition \ref{Prop2}, we have
\be\label{c-infty}
\lambda_\epsilon\ra +\infty,\quad
c_\epsilon\ra +\infty.
\ee
The convergence of $u_\epsilon$ will be described in the following proposition.

\begin{proposition}\label{Prop3} Assume $\alpha<\lambda_1(\Sigma)$ and $J_{\alpha,8\pi}$
has no minimizer in $\mathcal{H}$.
Let $u_\epsilon$
be a sequence of solutions to the equation (\ref{EL}). Let $c_\epsilon=u_\epsilon(x_\epsilon)$
be defined as in (\ref{maximum}) and assume that
$x_\epsilon\ra p\in\Sigma$.   If we define
\be\label{varphi}\varphi_\epsilon(y)={u}_\epsilon\le(\exp_{x_\epsilon}(r_\epsilon y)\ri)-c_\epsilon\ee
 and
\be\label{scal}r_\epsilon=\f{\sqrt{\lambda_\epsilon}}{\sqrt{8\pi(1-\epsilon)h(p)}}e^{-c_\epsilon/2},\ee
then
\be\label{bubble}\varphi_\epsilon(y)\ra\varphi(y)=-2\log(1+|y|^2/8)\quad{\rm in}\quad
C^1_{\rm loc}(\mathbb{R}^2).\ee
Moreover, $u_\epsilon$ converges to a Green function $G_p$ weakly in $W^{1,q}(\Sigma)$ for any $1<q<2$,
strongly in $L^r(\Sigma)$ for all $0<r<2q/(2-q)$, and in $C^1_{\rm loc}(\Sigma\setminus\{p\})$, where
$G_p$ satisfies
\be\label{Green}
\le\{\begin{array}{lll}
\Delta_{g}G_p(x)-\alpha G_p(x)=8\pi\delta_{p}(x)-\f{8\pi}{{\rm Vol}_g(\Sigma)}\\[1.5ex]
\int_\Sigma G_p(x)dv_{g}=0.
\end{array}
\ri.
\ee
\end{proposition}

The proof of Proposition \ref{Prop3} will be divided into several lemmas.
\begin{lemma}\label{r-0}
Let $r_\epsilon$ be defined as in (\ref{scal}).
For any $\gamma<1/2$, there holds $r_\epsilon^2 e^{\gamma c_\epsilon}\ra 0$.
In particular, $r_\epsilon c_\epsilon^q\ra 0$ for any $q>0$.
\end{lemma}
\noindent{\it Proof.} Multiplying both sides of the equation (\ref{EL}) by
$u_\epsilon$, we have
\bea\nonumber
\|u_\epsilon\|_{1,\alpha}^2&=&\int_\Sigma(|\nabla_gu_\epsilon|^2-\alpha u_\epsilon^2)dv_g\\
\nonumber
&=&\f{8\pi(1-\epsilon)}{\lambda_\epsilon}\int_\Sigma hu_\epsilon e^{u_\epsilon}dv_g\\ \label{energy}
&\leq&8\pi c_\epsilon.
\eea
In view of the Trudinger-Moser inequality (\ref{strong-TM}), we estimate
\bna
\int_\Sigma he^{u_\epsilon}dv_g&\leq&C\int_\Sigma e^{4\pi\f{u_\epsilon^2}{\|u_\epsilon\|_{1,\alpha}^2}+
\f{\|u_\epsilon\|_{1,\alpha}^2}{16\pi}}dv_g\\
&\leq&Ce^{\f{\|u_\epsilon\|_{1,\alpha}^2}{16\pi}}
\leq Ce^{\f{1}{2}c_\epsilon}.
\ena
It follows that
$$r_\epsilon^2=\f{\int_\Sigma he^{u_\epsilon}dv_g}{8\pi(1-\epsilon)h(p)}e^{-c_\epsilon}\leq Ce^{-\f{1}{2}c_\epsilon}.$$
This together with (\ref{c-infty}) gives the desired result. $\hfill\Box$\\

Let $0<\delta< i_g(\Sigma)$ be fixed and $i_g(\Sigma)$ be the injectivity radius of $(\Sigma,g)$. For
$y\in \mathbb{B}_{\delta r_\epsilon^{-1}}(0)$, the Euclidean ball of center $0$ and radius $\delta r_\epsilon^{-1}$, we set
\bea\label{psi}&&\psi_\epsilon(y)=c_\epsilon^{-1}{u}_\epsilon\le(\exp_{{x}_\epsilon}(r_\epsilon y)\ri),\\
&&g_\epsilon(y)=\le(\exp^\ast_{x_\epsilon}g\ri)(r_\epsilon y).\label{g-eps}\eea
Clearly $g_\epsilon\ra g_0$, the standard Euclidean metric, in $C^2_{\rm loc}(\mathbb{R}^2)$ as $\epsilon\ra 0$.
Note that $\psi_\epsilon\leq \psi_\epsilon(0)=1$. Concerning the asymptotic behavior of $\psi_\epsilon$, we have the following:
\begin{lemma}\label{psi-lemma}
$\psi_\epsilon\ra 1$ in $C^1_{\rm loc}(\mathbb{R}^2)$.
\end{lemma}
\noindent {\it Proof.} In view of (\ref{EL}), (\ref{psi}) and (\ref{g-eps}), we have
\bea\label{psi-equ}
\Delta_{g_\epsilon}\psi_\epsilon(y)
=\alpha r_\epsilon^2\psi_\epsilon(y)+c_\epsilon^{-1}\f{h(\exp_{{x}_\epsilon}(r_\epsilon y))}{h(p)}
e^{{u}_\epsilon(\exp_{{x}_\epsilon}(r_\epsilon y))-c_\epsilon}
-\f{8\pi(1-\epsilon)}{{\rm Vol}_g(\Sigma)}r_\epsilon^2 c_\epsilon^{-1}.
\eea
Let $q>1$ be any fixed number. By (\ref{energy}) and the Sobolev embedding theorem, we have
\be\label{Sobo}\int_\Sigma|u_\epsilon|^qdv_g\leq C\|u_\epsilon\|_{1,\alpha}^{q}\leq Cc_\epsilon^{\f{q}{2}}.\ee
Let $B_r(x)$ be a geodesic ball centered at $x\in \Sigma$ with radius $r$. It follows from a change of variables,
(\ref{Sobo}) and Lemma \ref{r-0}  that
\bea\nonumber
\int_{\mathbb{B}_R(0)}|r_\epsilon^2\psi_\epsilon(y)|^qdy&=&(1+o_\epsilon(1))\int_{B_{Rr_\epsilon}(x_\epsilon)}
r_\epsilon^{2q-2}c_\epsilon^{-q}|u_\epsilon|^qdv_g\\\nonumber
&\leq&Cr_\epsilon^{2q-2}c_\epsilon^{-q}\int_\Sigma|u_\epsilon|^qdv_g\\\nonumber
&\leq&Cr_\epsilon^{2q-2}c_\epsilon^{-q}\|u_\epsilon\|_{1,\alpha}^q\\\label{Lq}
&\leq&Cr_\epsilon^{2q-2}c_\epsilon^{-\f{q}{2}}\ra 0
\eea
as $\epsilon\ra 0$. Therefore we conclude that $\Delta_{g_\epsilon}\psi_\epsilon(y)$ converges to $0$ in $L^q_{\rm loc}(\mathbb{R}^2)$ for any
$q>1$. Noting that $\psi_\epsilon(y)\leq 1$ for all $y\in\mathbb{B}_{\delta r_\epsilon^{-1}}(0)$ and applying elliptic estimates to (\ref{psi-equ}),
we obtain $\psi_\epsilon\ra\psi$ in $C^1_{\rm loc}(\mathbb{R}^2)$ for some $\psi$ satisfying
$$\le\{\begin{array}{lll}
-\Delta_{\mathbb{R}^2}\psi=0\quad{\rm in}\quad\mathbb{R}^2\\[1.5ex]
\psi(y)\leq\psi(0)=1,
\end{array}\ri.$$
where $\Delta_{\mathbb{R}^2}$ denotes the usual Laplacian operator on $\mathbb{R}^2$. Then the Liouville theorem leads to
$\psi(y)\equiv 1$ for $y\in\mathbb{R}^2$. This completes the proof of the lemma. $\hfill\Box$\\

Let $\varphi_{\epsilon}(y)$ be defined as in (\ref{varphi})
for $y\in\mathbb{B}_{\delta r_\epsilon^{-1}}(0)$. To prove (\ref{bubble}), we calculate on $\mathbb{B}_{\delta r_\epsilon^{-1}}(0)$,
\bea\label{phi-equa}
\Delta_{g_\epsilon}\varphi_\epsilon(y)=\alpha r_\epsilon^2{u}_\epsilon\le(\exp_{x_\epsilon}(r_\epsilon y)\ri)
+\f{h(\exp_{x_\epsilon}(r_\epsilon y))}{h(p)}e^{\varphi_\epsilon(y)}
-\f{8\pi(1-\epsilon)}{{\rm Vol}_g(\Sigma)}r_\epsilon^2.
\eea
An obvious analog of (\ref{Lq}) implies that $r_\epsilon^2{u}_\epsilon(\exp_{x_\epsilon}(r_\epsilon y))\ra 0$
 in $L^q_{\rm loc}(\mathbb{R}^2)$  as $\epsilon\ra 0$ for any $q>1$.
Note that $\varphi_\epsilon(y)\leq \varphi_\epsilon(0)=0$. In view of Lemma \ref{psi-lemma},
we have by applying elliptic estimates to (\ref{phi-equa}),
$\varphi_\epsilon\ra \varphi$ in $C^1_{\rm loc}(\mathbb{R}^2)$ as $\epsilon\ra 0$, where $\varphi$ satisfies
$$\le\{\begin{array}{lll}
-\Delta_{\mathbb{R}^2}\varphi(y)=e^{\varphi(y)},\quad y\in\mathbb{R}^2\\[1.5ex]
\int_{\mathbb{R}^2}e^{\varphi(y)}dy<\infty.\end{array}
\ri.$$
A result of Chen-Li \cite{CL} implies that $\varphi$ can be written as in (\ref{bubble}) and thus
\be\label{bubble-energy}\int_{\mathbb{R}^2}e^{\varphi(y)}dy=8\pi.\ee

\begin{lemma}\label{delta} There holds
 $\lambda_\epsilon^{-1}he^{u_\epsilon}\rightharpoonup \delta_p$ in sense of measure, where $\delta_p$ denotes
the Dirac measure centered at $p$.
\end{lemma}

\noindent{\it Proof.} By a change of variables, we have
\bna
\f{8\pi(1-\epsilon)}{\lambda_\epsilon}\int_{B_{Rr_\epsilon}(x_\epsilon)}he^{u_\epsilon}dv_g
&=&(1+o_\epsilon(1))\int_{\mathbb{B}_R(0)}e^{\varphi_\epsilon(y)}dy\\
&=&(1+o_\epsilon(1))\int_{\mathbb{B}_R(0)}e^{\varphi(y)}dy.
\ena
This together with (\ref{bubble-energy}) leads to
\be\label{lim-1}\lim_{R\ra\infty}\lim_{\epsilon\ra 0}\int_{B_{Rr_\epsilon}(x_\epsilon)}\lambda_\epsilon^{-1} he^{u_\epsilon}dv_g=1.\ee
Hence
\be\label{lim-0}\lim_{R\ra\infty}\lim_{\epsilon\ra 0}\int_{\Sigma\setminus B_{Rr_\epsilon}(x_\epsilon)}
\lambda_\epsilon^{-1} he^{u_\epsilon}dv_g=0.\ee
Combining (\ref{lim-1}) and (\ref{lim-0}), we have for any $\eta\in C^0(\Sigma)$,
\bna
\lambda_\epsilon^{-1}\int_\Sigma \eta he^{u_\epsilon}dv_g
=\eta(p)+o_\epsilon(1).
\ena
This gives the desired result. $\hfill\Box$
\begin{lemma}\label{1q}
If $u\in C^2(\Sigma)$ is a solution of $\Delta_gu=f$, then for any $1<q<2$,
there exists some constant $C$ depending only on $(\Sigma,g)$ and $q$ such that
$$\|\nabla_gu\|_{L^q(\Sigma)}\leq C\|f\|_{L^1(\Sigma)}.$$
\end{lemma}

\noindent {\it Proof.} Without loss of generality we assume $\int_\Sigma udv_g=0$. Let $G(x,y)$ be the standard Green function
satisfying $\Delta_{g,y} G(x,y)=\delta_x(y)-1/V$ and $\int_\Sigma G(x,y) dv_{g,y}=0$, where $\delta_x(y)$ denotes the Dirac measure
centered at $x$, and $V$ is the area of $\Sigma$. Clearly
$$u(x)=\int_\Sigma f(y) G(x,y)dv_{g,y}.$$
It follows from \cite{Fontana} that $|\nabla_{g,x}G(x,y)|\leq C r(x,y)^{-1}$ for some constant $C$ depending only on $(\Sigma,g)$,
where $r(x,y)$ stands for the geodesic distance between $x$ and $y$.
Let $1<q<2$. One calculates by using the H\"older inequality
$$|\nabla_g u(x)|^q\leq C\|f\|_{L^1(\Sigma)}^{q-1}\int_\Sigma|f(y)|r(x,y)^{-q}dv_{g,y}.$$
Hence $\|\nabla_g u\|_{L^q(\Sigma)}\leq C\|f\|_{L^1(\Sigma)}$ for some constant $C$ depending only on $(\Sigma,g)$ and $q$.
$\hfill\Box$

\begin{lemma}\label{W1q}
For any $1<q<2$, there exists some constant $C$ such that $\|\nabla_gu_\epsilon\|_{L^q(\Sigma)}\leq C$.
\end{lemma}

\noindent{\it Proof.} Clearly (\ref{EL}) gives
\be\label{f-e}\Delta_gu_\epsilon=f_\epsilon=\alpha u_\epsilon+\f{8\pi(1-\epsilon)}
{\lambda_\epsilon}he^{u_\epsilon}-\f{8\pi(1-\epsilon)}{{\rm Vol}_g(\Sigma)}.\ee
In view of Lemma \ref{1q}, it suffices to prove
\be\label{L1}
\|f_\epsilon\|_{L^1(\Sigma)}\leq C.
\ee
Noting that $\lambda_\epsilon^{-1}\int_\Sigma he^{u_\epsilon}dv_g=1$, we only need to prove that $u_\epsilon$
is bounded in $L^1(\Sigma)$. For otherwise we can assume $\|u_\epsilon\|_{L^1(\Sigma)}\ra +\infty$ as $\epsilon\ra 0$.
Set $v_\epsilon=u_\epsilon/\|u_\epsilon\|_{L^1(\Sigma)}$. Then $\|v_\epsilon\|_{L^1(\Sigma)}=1$ and
$\Delta_gv_\epsilon=f_\epsilon/\|u_\epsilon\|_{L^1(\Sigma)}$. Obviously $f_\epsilon/\|u_\epsilon\|_{L^1(\Sigma)}$ is bounded in
$L^1(\Sigma)$. Given any $1<q<2$. It follows from Lemma \ref{1q} that $v_\epsilon$ is bounded in $W^{1,q}(\Sigma)$.
One can assume up to a subsequence, $v_\epsilon$ converges to $v$ weakly in $W^{1,q}(\Sigma)$, strongly in
$L^r(\Sigma)$ for any $0<r<2q/(2-q)$, and almost everywhere in $\Sigma$. Moreover, $v$ is a distributional solution of
$$\label{contrary-v}
\le\{\begin{array}{lll}\Delta_gv-\alpha v=0
\\[1.5ex]
\int_\Sigma v dv_g=0.\end{array}
\ri.
$$
This leads to $v\equiv 0$ contradicting the fact that
$\|v\|_{L^1(\Sigma)}=\lim_{\epsilon\ra 0}\|v_\epsilon\|_{L^1(\Sigma)}=1$.
Therefore $\|u_\epsilon\|_{L^1(\Sigma)}$ must be bounded and thus (\ref{L1}) holds. $\hfill\Box$\\

Combining Lemmas \ref{delta} and \ref{W1q}, we obtain for all $1<q<2$ and $0<r<2q/(2-q)$,
\bna
&u_\epsilon\rightharpoonup G_p&\quad{\rm weakly\,\,in}\quad W^{1,q}(\Sigma)\\
&u_\epsilon\rightarrow G_p&\quad{\rm strongly\,\,in}\quad L^r(\Sigma)\\
&u_\epsilon\rightarrow G_p&\quad{\rm a.\,e.\,\,in}\quad \Sigma,
\ena
where $G_p$ is a distributional solution of (\ref{Green}). Applying elliptic estimates to (\ref{Green}), we have that
$G_p$ takes the form
\be\label{local}G_p(x)=-4\log r+A_p+\psi(x),\ee
where $r$ denotes the geodesic distance between $x$ and $p$, $\psi\in C^1(\Sigma)$ and $\psi(p)=0$.
To complete the proof of Proposition \ref{Prop3}, we also need the following:
\begin{lemma}\label{smooth-Gr}
$u_\epsilon \ra G_p$ in $C^1_{\rm loc}(\Sigma\setminus\{p\})$ as $\epsilon\ra 0$.
\end{lemma}

\noindent{\it Proof.} For any domain $\Omega\subset\subset\Sigma\setminus\{p\}$, let $u_\epsilon^{(1)}$ be a solution of
$$\label{u(1)}
\le\{\begin{array}{lll}\Delta_gu_\epsilon^{(1)}=\f{8\pi(1-\epsilon)}{\lambda_\epsilon}he^{u_\epsilon}\quad&{\rm in}\quad&\Omega
\\[1.5ex]
u_\epsilon^{(1)}=0&{\rm on}&\p\Omega.\end{array}
\ri.
$$
By Lemma \ref{delta}, $\lambda_\epsilon^{-1}he^{u_\epsilon}$ converges to $0$ in $L^1(\Omega)$ as $\epsilon\ra 0$.
A result of Brezis-Merle \cite{BM91} implies that for any $r>0$, there exists some constant $C$ such that
\be\label{u1-b}\|e^{|u_\epsilon^{(1)}|}\|_{L^r(\Omega)}\leq C.\ee
Setting $u_\epsilon^{(2)}=u_\epsilon-u_\epsilon^{(1)}$, we have on $\Omega$,
\be\label{u2}\Delta_gu_\epsilon^{(2)}=\alpha u_\epsilon-\f{8\pi(1-\epsilon)}{{\rm Vol}_g(\Sigma)}.\ee
In view of Lemma \ref{W1q}, $u_\epsilon$ is bounded in $L^r(\Sigma)$ for any $r>0$. For any $\Omega^\prime\subset\subset\Omega$,
applying elliptic estimates to (\ref{u2}),
we have that $u_\epsilon^{(2)}$ is uniformly bounded in $\Omega^\prime$. This together with (\ref{u1-b}) leads to
$\|f_\epsilon\|_{L^r(\Omega^\prime)}\leq C$
for some $r>2$, where $f_\epsilon$ is defined as in (\ref{f-e}). Then we get the desired result by applying elliptic estimates to
(\ref{f-e}). $\hfill\Box$

\subsection{Lower bound estimate}\label{2.5}
In this subsection, we shall derive a lower bound of $J_{\alpha,8\pi}$ on $\mathcal{H}$. Let $G_{x_\epsilon}$
be the Green function satisfying $\int_\Sigma G_{x_\epsilon}dv_g=0$ and
\be\label{Green-usual}\Delta_gG_{x_\epsilon}-\alpha G_{x_\epsilon}=8\pi\delta_{x_\epsilon}-\f{8\pi}{{\rm Vol}_g(\Sigma)}.\ee
 Clearly $G_{x_\epsilon}$ can be represented by
 \be\label{Gr-x-e}G_{x_\epsilon}=-4\log r+A_{x_\epsilon}+O(r),\ee
 where $r$ denotes the geodesic distance between $x_\epsilon$ and $x$. Moreover, there holds
 $A_{x_\epsilon}\ra A_p$ as $\epsilon\ra 0$.
 Similar to \cite{DJLW}, we have the following maximum principle.
\begin{lemma}\label{max-prin} For any $\epsilon>0$ and $R>0$, there exists a constant $L_\epsilon$ such that for all
$y\in\Sigma\setminus B_{Rr_\epsilon}(x_\epsilon)$, there holds
$$u_\epsilon(y)-G_{x_\epsilon}(y)\geq L_\epsilon= -c_\epsilon+2\log\lambda_\epsilon
-2\log\pi-2\log h(p)-A_{x_\epsilon}+o_\epsilon(1)+o_R(1).$$
\end{lemma}
\noindent{\it Proof}. Note that
$$\Delta_g(u_\epsilon-G_{x_\epsilon})-\alpha(u_\epsilon-G_{x_\epsilon})\geq 0\quad{\rm in}\quad \Sigma\setminus B_{Rr_\epsilon}(x_\epsilon).$$
In view of Proposition \ref{Prop3} and the formula (\ref{Gr-x-e}), we have on $\p B_{Rr_\epsilon}(x_\epsilon)$ that
$$u_\epsilon-G_{x_\epsilon}=-c_\epsilon+2\log\lambda_\epsilon
-2\log\pi-2\log h(p)-A_{x_\epsilon}+o_\epsilon(1)+o_R(1).$$
The desired result follows from the maximum principle immediately. $\hfill\Box$\\

For any fixed $R>0$, we have
\begin{align*}
\int_{\S}|\na_g u_\e|^2dv_g=\int_{B_{Rr_\e}(x_\e)}|\na_g u_\e|^2dv_g+\int_{\S\setminus{B_{Rr_\e}(x_\e)}}|\na_g u_\e|^2dv_g.
\end{align*}
By Proposition \ref{Prop3}  we have
\begin{align}\label{in}
\int_{B_{Rr_\e}(x_\e)}|\na_g u_\e|^2dv_g
=&\int_{B_R(0)}|\na \varphi|^2dx+o_\e(1)\nonumber\\
=&16\pi\log\le(1+\f{R^2}{8}\ri)-16\pi+o_\e(1)+o_R(1).
\end{align}
It follows from (\ref{EL}) that
\begin{align}
\label{out}
\int_{\S\setminus{B_{Rr_\e}(x_\e)}}|\na_g u_\e|^2dv_g
=&\int_{\S\setminus{B_{Rr_\e}(x_\e)}}u_\e \la_g u_\e dv_g-\int_{\p B_{Rr_\e}(x_\e)}u_\e\f{\p u_\e}{\p n}ds_g\nonumber\\
=&\int_{\S\setminus{B_{Rr_\e}(x_\e)}}\a u_\e^2dv_g+8\pi(1-\e)\int_{\S\setminus{B_{Rr_\e}(x_\e)}}\f{u_\e he^{u_\e}}{\lam_\e}dv_g\nonumber\\
  &-\f{8\pi(1-\e)}{{\rm Vol}_g(\S)}\int_{\S\setminus{B_{Rr_\e}(x_\e)}}u_\e dv_g-\int_{\p B_{Rr_\e}(x_\e)}u_\e\f{\p u_\e}{\p n}ds_g.
\end{align}
Lemma \ref{max-prin} leads to
\begin{align}
\label{out0}
&8\pi(1-\e)\int_{\S\setminus{B_{Rr_\e}(x_\e)}}\f{u_\e he^{u_\e}}{\lam_\e}dv_g\nonumber\\
\geq &8\pi(1-\e)\int_{\S\setminus{B_{Rr_\e}(x_\e)}}G_{x_\e}\f{ he^{u_\e}}{\lam_\e}dv_g-8\pi(1-\e)\int_{\S\setminus{B_{Rr_\e}(x_\e)}}\le(c_\e-2\log\lam_\e\ri)\f{ he^{u_\e}}{\lam_\e}dv_g\nonumber\\
&+ 8\pi(1-\e)\int_{\S\setminus{B_{Rr_\e}(x_\e)}}\le(-2\log\pi-2\log h(p)-A_{x_\e}+o_\epsilon(1)+o_R(1)\ri)\f{ he^{u_\e}}{\lam_\e}dv_g.
\end{align}
We estimate three terms on the right hand side of (\ref{out0}) respectively. Using (\ref{EL}) and (\ref{Green-usual}), we obtain
\begin{align}
\label{out1}
&8\pi(1-\e)\int_{\S\setminus{B_{Rr_\e}(x_\e)}}G_{x_\e}\f{ he^{u_\e}}{\lam_\e}dv_g\nonumber
\\=&\int_{\S\setminus{B_{Rr_\e}(x_\e)}}G_{x_\e}\le(\la_gu_\e-\a u_\e+\f{8\pi(1-\e)}{{\rm Vol}_g(\S)}\ri)dv_g\nonumber\\
 =&\f{8\pi}{{\rm Vol}_g(\S)}\int_{B_{Rr_\e}(x_\e)}u_\e dv_g-\f{8\pi(1-\e)}{{\rm Vol}_g(\S)}\int_{B_{Rr_\e}(x_\e)}G_{x_\e}dv_g\nonumber\\
 &+\int_{\p B_{Rr_\e}(x_\e)}G_{x_\e}\f{\p u_\e}{\p n}ds_g-\int_{\p B_{Rr_\e}(x_\e)}u_\e \f{\p G_{x_\e}}{\p n}ds_g.
\end{align}
It also follows from  (\ref{EL}) that
\begin{align}
\label{out2}
&8\pi(1-\e)\int_{\S\setminus{B_{Rr_\e}(x_\e)}}\le(c_\e-2\log\lam_\e\ri)\f{ he^{u_\e}}{\lam_\e}dv_g\nonumber
\\=&(-c_\e+2\log\lam_\e)\int_{\p B_{Rr_\e}(x_\e)}\f{\p u_\e}{\p n}ds_g+(-c_\e+2\log\lam_\e)\a\int_{B_{Rr_\e}(x_\e)}u_\e dv_g\nonumber\\
    &+\f{8\pi(1-\e)}{{\rm Vol}_g(\S)}(-c_\e+2\log\lam_\e)\le({\rm Vol}_g(\S)-{\rm Vol}_g\le(B_{Rr_\e}(x_\e)\ri)\ri).
\end{align}
In view of (\ref{lim-0}), one has
\begin{align}
\label{out3}
 8\pi(1-\e)\int_{\S\setminus{B_{Rr_\e}(x_\e)}}\le(-2\log\pi-2\log h(p)-A_{x_\e}+o_\epsilon(1)+o_R(1)\ri)\f{ he^{u_\e}}{\lam_\e}dv_g=o_\e(1).
\end{align}
Here, in (\ref{out3}), we use the fact that $A_{x_\epsilon}$ is bounded. Inserting (\ref{out1})-(\ref{out3}) into (\ref{out0}), we get the lower bound estimate of
${\lam_\e^{-1}}\int_{\S\setminus{B_{Rr_\e}(x_\e)}}{u_\e he^{u_\e}}dv_g$. Then inserting this lower bound to
(\ref{out}), we obtain the estimate of $\int_{\S\setminus{B_{Rr_\e}(x_\e)}}|\na_g u_\e|^2dv_g$ as below.
\begin{align}
\label{out20}
\int_{\S\setminus{B_{Rr_\e}(x_\e)}}|\na_g u_\e|^2dv_g
\geq&\int_{\S\setminus{B_{Rr_\e}(x_\e)}}\a u_\e^2dv_g+\f{8\pi}{{\rm Vol}_g(\S)}\int_{B_{Rr_\e}(x_\e)}u_\e dv_g-\f{8\pi(1-\e)}{{\rm Vol}_g(\S)}\int_{B_{Rr_\e}(x_\e)}G_{x_\e}dv_g\nonumber\\
&+\int_{\p B_{Rr_\e}(x_\e)}G_{x_\e}\f{\p u_\e}{\p n}ds_g
      -\int_{\p B_{Rr_\e}(x_\e)}u_\e \f{\p G_{x_\e}}{\p n}ds_g\nonumber\\
      &+(-c_\e+2\log\lam_\e)\int_{\p B_{Rr_\e}(x_\e)}\f{\p u_\e}{\p n}ds_g+(-c_\e+2\log\lam_\e)\a\int_{B_{Rr_\e}(x_\e)}u_\e dv_g\nonumber\\
    &+\f{8\pi(1-\e)}{{\rm Vol}_g(\S)}(-c_\e+2\log\lam_\e)\le({\rm Vol}_g(\S)-{\rm Vol}_g\le(B_{Rr_\e}(x_\e)\ri)\ri)\nonumber\\
  &+\f{8\pi(1-\e)}{{\rm Vol}_g(\S)}\int_{B_{Rr_\e}(x_\e)}u_\e dv_g-\int_{\p B_{Rr_\e}(x_\e)}u_\e\f{\p u_\e}{\p n}ds_g+o_\e(1).
\end{align}
Using Proposition \ref{Prop3} and Lemma \ref{max-prin}, one has
\begin{align}
\label{out21}
&-\int_{\p B_{Rr_\e}(x_\e)}\f{\p u_\e}{\p n}\le(u_\e-G_{x_\e}+c_\e-2\log\lam_\e\ri)ds_g\nonumber\\
\geq& \f{\pi R^2}{1+{R^2}/{8}}\le(-2\log\pi-2\log h(p)-A_{x_\e}\ri)+o_\e(1)+o_R(1).
\end{align}
In view of (\ref{local}) and Proposition \ref{Prop3}, we have
\begin{align*}
\label{}
-\int_{\p B_{Rr_\e}(x_\e)}u_\e \f{\p G_{x_\e}}{\p n}ds_g
=-\le(c_\e-2\log\le(1+\f{R^2}{8}\ri)+o_\e(1)\ri)\le(-8\pi+O(Rr_\e)\ri).
\end{align*}
Then it follows by Lemma \ref{r-0} that
\begin{align}
\label{out22}
-\int_{\p B_{Rr_\e}(x_\e)}u_\e \f{\p G_{x_\e}}{\p n}ds_g
=8\pi c_\e-16\pi\log\le(1+\f{R^2}{8}\ri)+o_\e(1)+o_R(1).
\end{align}
Also Proposition \ref{Prop3} and Lemma \ref{r-0} lead to
\begin{align}
\label{out23}
\int_{B_{Rr_\e}(x_\e)}u_\e^2dv_g=o_\e(1),~~\int_{B_{Rr_\e}(x_\e)}u_\e dv_g=o_\e(1).
\end{align}
In view of (\ref{local}) and Lemma \ref{r-0}, we have
\begin{align}
\label{out24}
\int_{B_{Rr_\e}(x_\e)}G_{x_\e}dv_g=o_\e(1).
\end{align}
It follows from Lemma \ref{l>0}, Proposition \ref{Prop3} and Lemma \ref{r-0} that
\begin{align}
\label{out25}
(-c_\e+2\log\lam_\e)\int_{B_{Rr_\e}(x_\e)}u_\e dv_g=o_\e(1)
\end{align}
and
\begin{align}
\label{out26}
(-c_\e+2\log\lam_\e){\rm Vol}_g\le(B_{Rr_\e}(x_\e)\ri)=o_\e(1).
\end{align}
Hence we have by inserting (\ref{out21})-(\ref{out26}) into (\ref{out20}),
\begin{align}
\label{out100}
\int_{\S\setminus{B_{Rr_\e}(x_\e)}}|\na_g u_\e|^2dv_g
\geq& \a \int_\S u_\e^2dv_g+8\pi\le(-2\log\pi-2\log h(p)-A_{x_\e}\ri)+8\pi\e c_\e\nonumber\\
&+16\pi(1-\e)\log\lam_\e-16\pi\log\le(1+\f{R^2}{8}\ri)+o_\e(1)+o_R(1).
\end{align}
Combining (\ref{in}) and (\ref{out100}), we have
\begin{align*}
J_{\a, 8\pi(1-\e)}(u_\e)=&\f{1}{2}\int_{\S}|\na_g u_\e|^2dv_g-\f{\a}{2}\int_{\S}u_\e^2dv_g-8\pi(1-\e)\log\lam_\e\nonumber\\
\geq&-8\pi-8\pi\log\pi-8\pi\log h(p)-4\pi A_{x_\e}+8\pi\e c_\e+o_\e(1)+o_R(1).
\end{align*}
Letting $\e\ra0$ first and then $R\ra+\infty$, one has
\begin{align}
\label{lowerbound}
\inf_{u\in \mathcal{H}}J_{\alpha,8\pi}(u)
\geq &-8\pi-8\pi\log\pi-8\pi\log h(p)-4\pi A_{p}\nonumber\\
\geq &-8\pi-8\pi\log\pi-4\pi\max_{x\in\S}\le(2\log h(x)+A_{x}\ri).
\end{align}

\subsection{Test function computation}

In this subsection, we construct a sequence of functions $(\phi_\e)_{\e>0}$ satisfying
\be\label{74'}\lim_{\e\ra0}J_{\a, 8\pi}(\phi_\e-\bar{\phi}_\e)=-8\pi-8\pi\log\pi-4\pi\max_{x\in\S}\le(2\log h(x)+A_{x}\ri),\ee
where
$$\bar{\phi}_\e=\f{1}{{\rm Vol}_g(\Sigma)}\int_\Sigma \phi_\e dv_g.$$
Suppose that $2\log h(p)+A_p=\max_{x\in\S}\le(2\log h(x)+A_{x}\ri)$. Let $r=r(x)$ be the geodesic distance between $x$ and
$p$. We set
\be\label{p-eps}\phi_\epsilon(x)=\le\{
\begin{array}{lll}
c-2\log(1+\f{r^2}{8\epsilon^2}),&x\in B_{R\epsilon}(p)\\[1.5ex]
G_p(x)-\eta(x)\psi(x),&x\in B_{2R\epsilon}(p)\setminus B_{R\epsilon}(p)\\[1.5ex]
G_p(x),&x\in\Sigma\setminus B_{2R\epsilon}(p),
\end{array}\ri.\ee
where $\eta\in C_0^{\infty}(B_{2R\e}(p))$ is a cut-off function, $\eta\equiv 1$ in $B_{R{\e}}(p)$, $|\nabla_g\eta(x)|\leq\f{4}{R{\e}}$
for all $x\in B_{2R\e}(p)$, $\psi$ is defined as in (\ref{local}),
$$c=2\log(1+R^2/8)-4\log R-4\log\epsilon+A_p$$
and $R=R(\epsilon)$ satisfying $R\ra+\infty$ and $(R{\e})^2\log R\ra0$ as $\e\ra 0$.

A straightforward calculation shows
\begin{align}
\label{inner}
\int_{B_{R{\e}}(p)}|\na_g \phi_\e|^2dv_g=16\pi\log(1+R^2/8)-16\pi+o_\e(1).
\end{align}
Moreover we have
\begin{align}
\label{outer}
\int_{\Sigma\setminus B_{R{\e}}(p)}|\na_g \phi_\e|^2dv_g
=&\int_{\S\setminus B_{R{\e}}(p)}|\na_g G_p|^2dv_g+\int_{B_{2R{\e}}(p)\setminus B_{R{\e}}(p)}|\na_g (\eta\psi)|^2dv_g\nonumber\\
  &-2\int_{B_{2R{\e}}(p)\setminus B_{R{\e}}(p)}\na_g G_p \na_g(\eta\psi)dv_g\nonumber\\
=&-\int_{\p B_{R{\e}}(p)}G_p\f{\p G_p}{\p n}ds_g+\int_{\S\setminus B_{R{\e}}(p)}G_p\la_g G_pdv_g\nonumber\\
  &+\int_{B_{2R{\e}}(p)\setminus B_{R{\e}}(p)}|\na_g (\eta\psi)|^2dv_g\nonumber\\
  &-2\int_{B_{2R{\e}}(p)\setminus B_{R{\e}}(p)}\eta\psi\la_g G_p dv_g+2\int_{\p B_{R{\e}}(p)}\eta\psi\f{\p G_p}{\p n}ds_g.
\end{align}
Using (\ref{local}) one has
\begin{align}
\label{outer2}
-\int_{\p B_{R{\e}}(p)}G_p\f{\p G_p}{\p n}ds_g=-32\pi\log(R\e)+8\pi A_p+o_\e(1).
\end{align}
By (\ref{Green-usual}) and (\ref{local}) we have
\begin{align}
\label{outer4}
\int_{\S\setminus B_{R{\e}}(p)}G_p\la_g G_pdv_g
=&\int_{\S\setminus B_{R{\e}}(p)}G_p\le(\a G_p-\f{8\pi}{{\rm Vol}_g(\S)}\ri)dv_g\nonumber\\
=&\a \int_{\S} G^2_pdv_g+o_\e(1),
\end{align}
where we have used the facts $\int_{B_{R{\e}}(p)}G^2_pdv_g=o_\e(1)$ and $\int_{B_{R{\e}}(p)}G_pdv_g=o_\e(1)$.

Moreover one can easily get that
\begin{align}
\label{out6}
\int_{B_{2R{\e}}(p)\setminus B_{R{\e}}(p)}|\na_g (\eta\psi)|^2dv_g=o_\e(1),
\end{align}
that
\begin{align}
\label{outer8}
   &-2\int_{B_{2R{\e}}(p)\setminus B_{R{\e}}(p)}\eta\psi\la_g G_p dv_g
=-2\int_{B_{2R{\e}}(p)\setminus B_{R{\e}}(p)}\eta\psi\le(\a G_p-\f{8\pi}{{\rm Vol}_g(\S)}\ri)dv_g=o_\e(1),
\end{align}
and that
\begin{align}
\label{outer10}
2\int_{\p B_{R{\e}}(p)}\eta\psi\f{\p G_p}{\p n}ds_g=o_\e(1),
\end{align}
where we have used (\ref{Green-usual}) and (\ref{local}).
Inserting (\ref{outer2})-(\ref{outer10}) into (\ref{outer}) we obtain
\begin{align}
\label{outer100}
\int_{\S\setminus B_{R{\e}}(p)}|\na_g \phi_\e|^2dv_g=-32\pi\log(R \e)+8\pi A_p+\a\int_\S G^2_p dv_g+o_\e(1).
\end{align}
Combining (\ref{inner}) and (\ref{outer100}),  we have
\begin{align}
\label{energylast}
\int_\S |\na_g \phi_\e|^2dv_g=-32\pi\log\e-16\pi\log 8-16\pi+8\pi A_p+\a\int_\S G^2_p dv_g+o_\e(1).
\end{align}
Clearly
\begin{align}
\label{meanvalue}
\bar{\phi}_\e=\f{1}{{\rm Vol}_g(\S)}\int_{\S} \phi_\e dv_g=o_\e(1)
\end{align}
and
\begin{align}
\label{L2}
\int_\S \le(\phi_\e-\bar{\phi}_\e\ri)^2dv_g=\int_\S G^2_pdv_g+o_\e(1).
\end{align}

We now estimate $\int_\S he^{\phi_\e}dv_g$. Choosing $\d >0$ sufficiently small and noting that $G_p$ has the expression (\ref{local}) in $B_\d (p)$,
we have
\begin{align}
\label{hint}
\int_\S he^{\phi_\e}dv_g
=&h(p)\int_{B_{R{\e}}(p)} e^{\phi_\e}dv_g+\int_{B_{R{\e}}(p)} (h-h(p))e^{\phi_\e}dv_g\nonumber\\
  &+\int_{B_\d(p)\setminus B_{R{\e}}(p)} he^{\phi_\e}dv_g+\int_{\S\setminus B_\d(p)} he^{\phi_\e}dv_g.
\end{align}
A straightforward calculation gives
\begin{align}
\label{hint1}
h(p)\int_{B_{R{\e}}(p)} e^{\phi_\e}dv_g=8\pi h(p)e^{-2\log 8-2\log\epsilon+A_p+o_\epsilon(1)}
\end{align}
and
\begin{align}
\label{hint2}
 \int_{B_{R{\e}}(p)} (h-h(p))e^{\phi_\e}dv_g=o_\e(1)\,\epsilon^{-2}.
\end{align}
Also one has
\begin{align}
\label{hint3}
0<\int_{B_\d(p)\setminus B_{R{\e}}(p)} he^{\phi_\e}dv_g
\leq& C(\max_{\S}h)\int_{B_\d(p)\setminus B_{R{\e}}(p)} e^{G_p}dv_g\nonumber\\
\leq& C(\max_{\S}h)\le(-\f{1}{\d^2}+\f{1}{(R\e)^2}\ri)
\end{align}
and
\begin{align}
\label{hint4}
\int_{\S\setminus B_\d(p)} he^{\phi_\e}dv_g\leq C(\max_{\S}h)\int_{\S\setminus B_\d(p)} e^{G_p}dv_g.
\end{align}
Inserting (\ref{hint1})-(\ref{hint4}) into (\ref{hint}), we have
\bna
\int_\Sigma he^{\phi_\epsilon}dv_g&=&(1+o_\epsilon(1))8\pi h(p)e^{-2\log 8-2\log\epsilon+A_p}+C\\
&=&(1+o_\epsilon(1))8\pi h(p)e^{-2\log 8-2\log\epsilon+A_p},
\ena
where $C$ is a constant depending only on $(\Sigma,g)$, $\max_\Sigma h$ and $\delta$. Hence
\begin{align}
\label{hint100}
\log\int_\S he^{\phi_\e}dv_g=-\log 8+\log(\pi h(p))-2\log\epsilon+A_p+o_\epsilon(1).
\end{align}
Combining (\ref{energylast}), (\ref{meanvalue}), (\ref{L2}) and (\ref{hint100}), we have
\begin{align}
\label{limit-1}
J_{\a, 8\pi}(\phi_\e-\bar{\phi}_\e)
=&\f{1}{2}\int_\S \le(|\na_g \phi_\e|^2-\a(\phi_\e-\bar{\phi}_\e)^2\ri)dv_g-8\pi\log\int_\S he^{\phi_\e-\bar{\phi}_\e}dv_g\nonumber\\
=&-8\pi-8\pi\log\pi-4\pi\le(2\log h(p)+A_p\ri)+o_\e(1).
\end{align}
This implies (\ref{74'}), which together with (\ref{lowerbound}) completes the proof of $(i)$ of Theorem \ref{main-theorem}. $\hfill\Box$
\section{Proof of Theorem \ref{main-theorem-2}}\label{sec2}

In this section we prove Theorem \ref{main-theorem-2} by using similar method of the proof of Theorem \ref{main-theorem}. Let
$\{e_i\}_{i=1}^{m_\ell}\subset\mathcal{H}\cap C^1(\Sigma)$ be
an orthonormal basis of $E_\ell=E_{\lambda_1(\Sigma)}\oplus \cdots\oplus E_{\lambda_\ell(\Sigma)}$, namely $E_\ell={\rm span}
\{e_1,\cdots,e_{m_\ell}\}$ and
$$\langle e_i,e_j\rangle=\int_\Sigma e_ie_jdv_g=\delta_{ij}=\le\{\begin{array}{lll}
1,&i=j\\[1.5ex]
0,&i\not=j
\end{array}\ri.$$
for all $i,j=1,\cdots,m_\ell$. Note that $\mathcal{H}=E_\ell\oplus E_\ell^\perp$. $v\in E_\ell^\perp$ if and only if
$v\in \mathcal{H}$ and $\langle v,e_i\rangle=0$ for all $i=1,\cdots,m_\ell$.\\

{\it Proof of $(i)$ of Theorem \ref{main-theorem-2}}. Let $\alpha<\lambda_{\ell+1}$ be fixed. Using the argument of Subsections \ref{2.3}-\ref{2.5} with minor modifications,
we obtain an analog of (\ref{lowerbound}) as the following: If $J_{\alpha,8\pi}$ has no minimizer in $E_\ell^\perp$, then
\be\label{lower-b3}\inf_{u\in E_{\ell}^\perp}J_{\alpha,8\pi}(u)
\geq -8\pi-8\pi\log\pi-4\pi\max_{x\in\S}\le(2\log h(x)+A_{\alpha,x}\ri),
\ee
where $A_{\alpha,x}=\lim_{y\ra x}(G_{\alpha,x}(y)+4\log r)$ is a constant, $r$ denotes the geodesic distance between $y$ and $x$,
$G_{\alpha,x}$ is a Green function
satisfying (\ref{green-alpha}). Assume
$$2\log h(p)+A_{\alpha,p}=\max_{x\in\S}\le(2\log h(x)+A_{\alpha,x}\ri).$$
It follows from elliptic estimates that $G_{\alpha,p}$ can be written as
$$G_{\alpha,p}(x)=-4\log r+A_{\alpha,p}+\psi_{\alpha}(x),$$
where $\psi_\alpha\in C^1(\Sigma)$ and $r$ denotes the geodesic distance between $x$ and $p$.
We now prove that if $J_{\alpha,8\pi}$ has no minimizer in $E_\ell^\perp$, then
\be\label{lower-equal}\inf_{u\in E_{\ell}^\perp}J_{\alpha,8\pi}(u)
= -8\pi-8\pi\log\pi-4\pi\le(2\log h(p)+A_{\alpha,p}\ri).
\ee
 Similar to (\ref{p-eps}), we set
$$\phi_{\alpha,\epsilon}(x)=\le\{
\begin{array}{lll}
c-2\log(1+\f{r^2}{8\epsilon^2}),&x\in B_{R\epsilon}(p)\\[1.5ex]
G_{\alpha,p}(x)-\eta(x)\psi_\alpha(x),&x\in B_{2R\epsilon}(p)\setminus B_{R\epsilon}(p)\\[1.5ex]
G_{\alpha,p}(x),&x\in\Sigma\setminus B_{2R\epsilon}(p),
\end{array}\ri.$$
where $\eta\in C_0^{\infty}(B_{2R\e}(p))$ is a cut-off function, $\eta\equiv 1$ in $B_{R{\e}}(p)$, $|\nabla_g\eta(x)|\leq\f{4}{R{\e}}$
for all $x\in B_{2R\e}(p)$,
$$c=2\log(1+R^2/8)-4\log R-4\log\epsilon+A_{\alpha,p}$$
and $R=R(\epsilon)$ satisfying $R\ra+\infty$ and $(R{\e})^2\log R\ra0$ as $\e\ra 0$. Similar to (\ref{limit-1}), we derive
\be\label{limit-2}
J_{\alpha,8\pi}(\phi_{\alpha,\epsilon}-\bar{\phi}_{\alpha,\epsilon})=-8\pi-8\pi\log\pi-4\pi\le(2\log h(p)+A_{\alpha,p}\ri)+o_\e(1),
\ee
where
\be\label{m-2}\bar{\phi}_{\alpha,\epsilon}=\f{1}{{\rm Vol}_g(\Sigma)}\int_\Sigma \phi_{\alpha,\epsilon}dv_g=o_\epsilon(1).\ee
Define a new sequence of functions
$$\psi_{\alpha,\epsilon}=\phi_{\alpha,\epsilon}-\bar{\phi}_{\alpha,\epsilon}-\sum_{k=1}^{m_\ell}\langle
\phi_{\alpha,\epsilon}-\bar{\phi}_{\alpha,\epsilon},e_k\rangle e_k,$$
where $\{e_i\}_{i=1}^{m_\ell}$ is an orthonormal basis on $E_\ell\subset\mathcal{H}$. Clearly $\psi_{\alpha,\epsilon}\in E_\ell^\perp$.

Noting that
\bna
&&\int_{\Sigma\setminus B_{R\epsilon}(p)}G_{\alpha,p}e_k dv_g=-\int_{B_{R\epsilon}(p)}G_{\alpha,p}e_k dv_g=o_\epsilon(1),\\
&&\int_{B_{R\epsilon}(p)}\phi_{\alpha,\epsilon}e_k dv_g=o_\epsilon(1),\,\,
\int_{B_{2R\epsilon}(p)\setminus B_{R\epsilon}(p)}\eta \psi_{\alpha}e_k dv_g=o_\epsilon(1),
\ena
we have
$$\langle\phi_{\alpha,\epsilon},e_k\rangle=\int_\Sigma \phi_{\alpha,\epsilon}e_k dv_g=o_\epsilon(1)$$
for any $k=1,\cdots,m_\ell$. This together with (\ref{m-2}) leads to
\be\label{tend-0}\langle\phi_{\alpha,\epsilon}-\bar{\phi}_{\alpha,\epsilon},e_k\rangle=
\int_\Sigma \phi_{\alpha,\epsilon}e_k dv_g+\bar{\phi}_{\alpha,\epsilon}\int_\Sigma e_k dv_g=o_\epsilon(1).\ee

Also we calculate
\bna
\int_{B_{R\epsilon}(p)}|\nabla_g\phi_{\alpha,\epsilon}|dv_g&=&(1+O((R\epsilon)^2))\int_0^{R\epsilon}\f{8\pi r^2}{r^2+8\epsilon^2}dr\\
&\leq& 8\pi(1+O((R\epsilon)^2))R\epsilon=o_\epsilon(1)
\ena
and
\bna
\int_{\Sigma\setminus B_{R\epsilon}(p)}|\nabla_g\phi_{\alpha,\epsilon}|dv_g&\leq&\int_{\Sigma}
|\nabla_gG_{\alpha,p}|dv_g+\int_{B_{2R\epsilon}(p)\setminus B_{R\epsilon}(p)}|\nabla_g(\eta\psi_\alpha)|dv_g\\
&\leq& C
\ena
for some constant $C$ depending only on $(\Sigma,g)$ and $G_{\alpha,p}$.
Hence we have
\be\label{grad-1}\int_\Sigma |\nabla_g\phi_{\alpha,\epsilon}|dv_g\leq C.\ee
Denote $\phi_{\alpha,\epsilon}^{(\ell)}=\sum_{k=1}^{m_\ell}\langle \phi_{\alpha,\epsilon}-
\bar{\phi}_{\alpha,\epsilon},e_k\rangle e_k$. Noting that $e_k\in C^1(\Sigma)$ for $k=1,\cdots,m_\ell$, we have by combining
(\ref{tend-0}) and (\ref{grad-1}) that
\bea\nonumber
\le|\int_\Sigma \nabla_g\phi_{\alpha,\epsilon}^{(\ell)}\nabla_g\phi_{\alpha,\epsilon}dv_g\ri|&\leq&
\int_\Sigma|\nabla_g\phi_{\alpha,\epsilon}|dv_g\sup_\Sigma|\nabla_g\phi_{\alpha,\epsilon}^{(\ell)}|\\\nonumber
&\leq& C\sum_{k=1}^{m_\ell}|\langle \phi_{\alpha,\epsilon}-
\bar{\phi}_{\alpha,\epsilon},e_k\rangle|\sup_\Sigma|\nabla_ge_k|\\\label{2-0}
&=&o_\epsilon(1).
\eea
Obviously (\ref{tend-0}) leads to
\be\label{t-0-2}\int_\Sigma|\nabla_g\phi_{\alpha,\epsilon}^{(\ell)}|^2dv_g=o_\epsilon(1).\ee
In view of (\ref{2-0}) and (\ref{t-0-2}), we have
\bea
\nonumber \int_\Sigma|\nabla_g\psi_{\alpha,\epsilon}|^2dv_g&=&\int_\Sigma |\nabla_g(\phi_{\alpha,\epsilon}-
\phi_{\alpha,\epsilon}^{(\ell)})|^2dv_g\\\nonumber
&=&\int_\Sigma|\nabla_g\phi_{\alpha,\epsilon}|^2dv_g+\int_\Sigma|\nabla_g\phi_{\alpha,\epsilon}^{(\ell)}|^2dv_g
-2\int_\Sigma\nabla_g\phi_{\alpha,\epsilon}\nabla_g\phi_{\alpha,\epsilon}^{(\ell)}dv_g\\\label{t-0-3}
&=&\int_\Sigma|\nabla_g\phi_{\alpha,\epsilon}|^2dv_g+o_\epsilon(1).
\eea
Noting that
$$\int_\Sigma|\phi_{\alpha,\epsilon}|dv_g\leq C\int_\Sigma |G_{\alpha,p}|dv_g$$
and recalling (\ref{tend-0}), we get
\bea\nonumber
\alpha\int_\Sigma\psi_{\alpha,\epsilon}^2dv_g&=&\alpha\int_\Sigma(\phi_{\alpha,\epsilon}-\bar{\phi}_{\alpha,\epsilon}
-\phi_{\alpha,\epsilon}^{(\ell)})^2dv_g\\\label{4-0}
&=&\alpha\int_\Sigma(\phi_{\alpha,\epsilon}-\bar{\phi}_{\alpha,\epsilon})^2dv_g+o_\epsilon(1).
\eea
Similarly we have
\be\label{5-0}
\log\int_\Sigma he^{\psi_{\alpha,\epsilon}}dv_g=\log\int_\Sigma he^{\phi_{\alpha,\epsilon}-\bar{\phi}_{\alpha,\epsilon}}dv_g+o_\epsilon(1).
\ee
It follows from (\ref{t-0-3})-(\ref{5-0}) that
$$J_{\alpha,8\pi}(\psi_{\alpha,\epsilon})=J_{\alpha,8\pi}(\phi_{\alpha,\epsilon}-\bar{\phi}_{\alpha,\epsilon})+o_\epsilon(1).$$
This together with (\ref{limit-2}) leads to
\be\label{limit-3}
J_{\alpha,8\pi}(\psi_{\alpha,\epsilon})=-8\pi-8\pi\log\pi-4\pi\le(2\log h(p)+A_{\alpha,p}\ri)+o_\e(1).
\ee
Comparing (\ref{limit-3}) with (\ref{lower-b3}), we conclude (\ref{lower-equal}) under the assumption that $J_{\alpha,8\pi}$
has no minimizer on $E_\ell^\perp$. This completes the proof of $(i)$ of Theorem \ref{main-theorem-2}. $\hfill\Box$\\

{\it Proof of $(ii)$ of Theorem \ref{main-theorem-2}}. Since the proof is completely analogous to
that of $(ii)$ of Theorem \ref{main-theorem},
we omit the details but leave it to the interested reader. $\hfill\Box$\\

{\it Proof of $(iii)$ of Theorem \ref{main-theorem-2}}. Let $\widetilde{M}_k\in \mathcal{H}$ be defined as in (\ref{M-tilde}).
We set
$$M_{\ell,k}(x)=\widetilde{M}_k(x)-\sum_{k=1}^{m_\ell}\langle \widetilde{M}_k,e_j\rangle e_j.$$
Obviously $M_{\ell,k}\in E_\ell^\perp$. One can check that
$$\log\int_\Sigma he^{M_{\ell,k}}dv_g=\le(\f{1}{2}+o_k(1)\ri)\log k$$
and
\bna
\int_\Sigma|\nabla_gM_{\ell,k}|^2dv_g-\alpha\int_\Sigma M_{\ell,k}^2dv_g&=&(1+O(r))4\pi\log k+O(1)\\
&\leq&(1+C_1r)4\pi\log k+C_2
\ena
for some positive constants $C_1$ and $C_2$ depending only on $(\Sigma,g)$, $\alpha$ and $r$. For any $\beta>8\pi$, we take $r>0$
such that $8\pi(1+C_1r)<\beta$. Then we have
$J_{\alpha,\beta}(M_{\ell,k})\ra -\infty$ as $k\ra+\infty$. This gives the desired result. $\hfill\Box$

\begin{rmk} 
  When we were students, we learned from Professor Jiayu Li  that they mistyped $x_\e$ as $p$ at some places in \cite{DJLW}.
  The first place is in Lemma 2.5, one needs choose a local normal coordinate system around $x_\e$. 
  The second place is Lemma 2.9,  it should be written as: In $\S\setminus B_{r_\e}(x_\e)$, we have 
  $u_\e\geq G_{x_\e}-c_\e-2\log(\f{1+\pi h(p)R^2}{R^2})-A_{x_\e}+o_\e(1)$, where $o_\e(1)\ra0$ as 
  $\e\ra0$. And the third place is in the proof of Lemma 2.2, the integral $\int_\S |\na u_\e|^2dv_g$ should be divided into
  $\int_{\S\setminus B_\d(x_\e)}|\na u_\e|^2dv_g+\int_{B_\d(x_\e)\setminus B_{r_\e}(x_\e)}|\na u_\e|^2dv_g
  +\int_{B_{r_\e}(x_\e)}|\na u_\e|^2dv_g$ and then be estimated respectively.
\end{rmk}

 \bigskip

  {\bf Acknowledgements}. Y. Yang is supported by the National Science Foundation of China
  (Grant Nos.11171347 and 11471014). X. Zhu is supported by the National Science Foundation of China
  (Grant Nos. 41275063 and 11401575).

\end{document}